\documentclass[12pt]{article}
\usepackage{graphicx,graphics,amssymb,amsfonts,amsmath,bm,url,color}
\usepackage[margin=2.5cm]{geometry}
\usepackage[english]{babel}

\def\no{\noindent}
\newtheorem{defi}{Definition}
\newtheorem{rem}{Remark}
\newtheorem{theo}{Theorem}

\def\proof{\underline{Proof.}~}
\def\QED{\mbox{$\Box{~}$}}
\def\pmatrix{\left(\begin{array}}
\def\endpmatrix{\end{array}\right)}
\def\RR{{\mathbb{R}}}
\def\CC{{\mathbb{C}}}
\def\DD{{\mathbb{D}}}
\def\bfb{{\bm{b}}}
\def\bfc{{\bm{c}}}
\def\bfe{{\bm{e}}}

\def\bfy{{\bm{y}}}
\def\bfgamma{{\bm{\gamma}}}
\def\hg{\hat{\gamma}}
\def\hbg{\hat{\bfgamma}}
\def\bfeta{\bm{\eta}}
\def\diag{{\rm diag}}
\def\PP{{\cal P}}
\def\XX{\hat{X}}
\def\aa{{\alpha}}

\def\dd{{\rm d}}
\def\hDelta{\hat{\Delta}}
\def\hc{\hat{c}}

\def\hPP{\hat{\PP}}

\def\hL{\hat{L}}
\def\hU{\hat{U}}
\def\hM{\hat{M}}
\def\hA{\hat{A}}

\title{Efficient implementation of Gauss collocation and Hamiltonian Boundary Value Methods}

\author{Luigi Brugnano\footnote{Dipartimento di Matematica e Informatica ``U.\,Dini'', Universit\`a di Firenze, Italy
({\tt luigi.brugnano@unifi.it})}
\and Gianluca Frasca Caccia\footnote{Dipartimento di Matematica e Informatica ``U.\,Dini'', Universit\`a di Firenze, Italy
({\tt frasca@math.unifi.it})}
\and Felice Iavernaro\footnote{Dipartimento di Matematica, Universit\`a di Bari, Italy
({\tt felice.iavernaro@uniba.it})}
}

\date{\small\em Warmly dedicated to celebrate the $80^{th}$ birthday of John Butcher}

\begin{document}
\maketitle

\begin{abstract}
In this paper we define an efficient implementation for the family of low-rank energy-conserving Runge-Kutta methods named Hamiltonian Boundary Value Methods (HBVMs), recently defined in the last years. The proposed implementation relies on the particular structure of the Butcher matrix defining such methods, for which we can derive an efficient splitting procedure. The very same procedure turns out to be automatically suited for the efficient implementation of Gauss-Legendre collocation methods, since these methods are a special instance of HBVMs. The linear convergence analysis of the splitting procedure exhibits excellent properties, which are confirmed by a few numerical tests.

\bigskip
\no{\bf Keywords:} Energy-conserving methods; Hamiltonian Boundary Value Methods; W-transform; Gauss-Legendre collocation methods; Implicit Runge-Kutta methods; Splitting.

\bigskip
\no{\bf MSC (2010):} 65P10, 65L05, 65L06, 65L99.
\end{abstract}

\section{Introduction}
The efficient numerical solution of implicit Runge-Kutta methods has been the subject of many investigations in the last
 decades, starting from the seminal papers of Butcher \cite{Bu76,Bu79} (see also \cite{CoBu83}). This aspect is even more relevant when dealing with {\em geometric} Runge-Kutta methods, that is, methods used in the framework of {\em Geometric Integration} where, usually, the discrete problems generated by the methods need to be solved to within full machine accuracy, in order not to waste the specific properties of the methods. 
 
 In more details, in this paper we shall deal with the numerical solution of Hamiltonian problems, namely problems in the form,
\begin{equation}\label{ivp}
y' = J\nabla H(y), \qquad y(t_0)=y_0\in\RR^{2m}, 
\end{equation}
where 
\begin{equation}\label{yJ}
y = \pmatrix{c} q\\ p\endpmatrix,\quad q,p\in\RR^m, \qquad J = \pmatrix{cc} O &I_m\\ -I_m &O\endpmatrix,
\end{equation}
$H(y)$ is  the  (scalar) {\em Hamiltonian} function defining the problem, and $I_m$ the identity matrix of dimension $m$.\footnote{In the following, when the size of the identity is not specified, it can be deduced from the context.} Due to the skew-symmetry of  matrix  $J$ one has
$$\frac{\dd}{\dd t}H(y(t)) = \nabla H(y(t))^Ty'(t) = \nabla H(y(t))^TJ\nabla H(y(t)) = 0,$$
so that $$H(y(t))=H(y_0), \qquad \forall\,t\ge t_0.$$ For isolated mechanical systems, the Hamiltonian has the physical meaning of the total energy of the system, so that often the Hamiltonian is referred to as the {\em energy}. Its conservation is, therefore, a significant feature for the discrete dynamical system induced by a numerical method for solving (\ref{ivp}): methods having this property are usually called {\em energy-conserving methods}. Among such methods, we are interested in the class of energy-conserving methods named {\em Hamiltonian Boundary Value Methods (HBVMs)} \cite{BIT10,BIT12_1} (see also \cite{BIT11,BIT12}, and  \cite{BI12,BIT12_3} for generalizations), which have been recently devised  starting from the concept of {\em discrete line integrals}, defined in \cite{IP07,IP08,IT09}.
For such methods, the discrete problem can be conveniently posed in a suitable form which can be exploited to derive efficient implementation strategies, as was done in \cite{BIT11}. Here we further improve on such results, by proposing and analysing an iterative procedure based on the particular structure of the discrete problem. As a by product, an efficient implementation of Gauss-Legendre collocation methods is also obtained. Indeed, these latter methods may be interpreted as a particular instance of HBVMs. The proposed procedure is strictly related to that recently devised in \cite{BIM12} for Radau IIA collocation formulae, though the two approaches are substantially different.  

With this premise, the paper is organized as follows: in Section~\ref{two} we describe the structure of the discrete problem generated by HBVMs, along with the way of solving it, as done so far; in Section~\ref{three} we introduce the new iterative procedure, which is based on a suitable splitting; in Section~\ref{four} we study the convergence properties of the new iteration, also comparing it with known existing ones; in Section~\ref{five} a few numerical tests are reported; at last, a few conclusions are contained in Section~\ref{six}.

\section{Discrete problem induced by HBVMs}\label{two}

We now recall the basic facts about HBVMs, and derive the most efficient formulation of the generated discrete problems. Let us assume, for sake of brevity, $t_0=0$ in (\ref{ivp}), and consider the approximation of the problem over the interval $[0,h]$, which will clearly concern the very first application of a given numerical method.
Let us then consider the orthonormal polynomial basis, on the interval $[0,1]$, provided by the shifted and scaled Legendre polynomials $\{P_j\}$:
\begin{equation}\label{orto}
\deg P_i = i, \qquad \int_0^1 P_i(x)P_j(x)\dd x = \delta_{ij}, \qquad \forall\, i,j\ge 0,
\end{equation}
where $\delta_{ij}$ is the usual Kronecker symbol. Under suitable mild assumptions on the Hamiltonian function $H$, the right-hand side of the differential equation (\ref{ivp}) can be expanded along the considered basis, thus giving
\begin{equation}\label{ylocal}
y'(ch) = \sum_{j\ge0} \gamma_j(y)P_j(c), \qquad c\in[0,1],
\end{equation}
where
\begin{equation}\label{gammaj}
\gamma_j(y) = \int_0^1 J\nabla H(y(\tau h))P_j(\tau)\dd\tau, \qquad j\ge0.
\end{equation} 
By imposing the initial condition, the solution of this problem is formally obtained by
\begin{equation}\label{ysol}
y(ch) = y_0 +h\sum_{j\ge0} \gamma_j(y)\int_0^c P_j(x)\dd x, \qquad c\in[0,1].\end{equation}
In order to derive a polynomial approximation $\sigma$ of degree $s$ to (\ref{ysol}), we consider the following approximated ODE-IVPs:
\begin{equation}\label{sigma1}
\sigma'(ch) = \sum_{j=0}^{s-1} \gamma_j(\sigma)P_j(c), \qquad c\in[0,1], \qquad \sigma(0)=y_0,
\end{equation}
where $\gamma_j(\sigma)$ is defined according to (\ref{gammaj}), by formally replacing $y$ by $\sigma$. Consequently, the approximation to (\ref{ysol}) will be given by
\begin{equation}\label{sigma}
\sigma(ch) = y_0 +h\sum_{j=0}^{s-1} \gamma_j(\sigma)\int_0^c P_j(x)\dd x, \qquad c\in[0,1].\end{equation}
 For sake of simplicity, assume now that the Hamiltonian function is a polynomial of degree $\nu$  (for the general case, see Theorem~\ref{propHBVM} below). Consequently, the (unknown) vector coefficients $\{\gamma_j(\sigma)\}$ can be exactly obtained by using a quadrature formula defined at the Gaussian abscissae $\{c_1,\dots,c_k\}$, i.e.,
\begin{equation}\label{ci}
P_k(c_i) = 0, \qquad i=1,\dots,k,
\end{equation}
and corresponding weights $\{b_1,\dots,b_k\}$,\footnote{Hereafter, we shall always assume this choice.}
\begin{equation}\label{quad}
\gamma_j(\sigma) = \sum_{i=1}^k b_iJ\nabla H(\sigma(c_i h))P_j(c_i), \qquad j=0,\dots,s-1,
\end{equation} provided that 
\begin{equation}\label{nu}
\nu\le\frac{2k}s.
\end{equation}
By setting 
\begin{equation}\label{Yi}
Y_i = \sigma(c_i h), \qquad i=1,\dots,k,
\end{equation}
and considering that the new approximation is given by
$$y_1\equiv \sigma(h) = y_0+h\int_0^1 J\nabla H(\sigma(\tau h))\dd\tau 
= y_0+h\sum_{i=1}^k b_i J\nabla H(Y_i),$$
one then obtains the following $k$-stage Runge-Kutta method,
\begin{equation}\label{RK}
\begin{array}{c|c} \bfc & A\\ \hline&\bfb^T\end{array}
\end{equation}
where, as usual, $\bfb,\bfc\in\RR^k$ are the vectors containing the weights and the abscissae, respectively, and (see, e.g. \cite{BIT11,BIT12,BIT12_1})
\begin{equation}\label{A}
A = \PP_{s+1}\XX_s\PP_s^T\Omega\in\RR^{k\times k},
\end{equation}
with  
\begin{eqnarray}\label{Pr}
\PP_r &=& \pmatrix{ccc}
P_0(c_1) & \dots & P_{r-1}(c_1)\\
\vdots   &           &\vdots\\
P_0(c_k) & \dots & P_{r-1}(c_k)\endpmatrix
\in\RR^{k\times r}, \qquad r=s,s+1,\\ \label{Xs}
\XX_s &=& \pmatrix{cccc}
\frac{1}2 & -\xi_1\\
\xi_1      &0          &\ddots\\
              &\ddots  &\ddots &-\xi_{s-1}\\
              &            &\xi_{s-1} &0 \\ \hline
              &            &               &\xi_{s} \endpmatrix~\equiv~\pmatrix{c} X_s\\ \hline
              0\dots 0\,\xi_s\endpmatrix\in\RR^{s+1\times s},\\
              \xi_i &=& \left(2\sqrt{4i^2-1}\right)^{-1}, \qquad i=1,\dots,s,\label{xi}\\
              \Omega &=& \diag(\bfb)\in\RR^{k\times k}.\label{Ome}
\end{eqnarray} 
We observe that, when $k=s$,  (\ref{A}) becomes the $W$-transformation \cite[pag.\,79]{HW96} of the $s$-stage Gauss-Legendre Runge-Kutta method. Consequently, (\ref{A}) can be also regarded as a generalization of the $W$-transformation.

Clearly, the Runge-Kutta method (\ref{RK})--(\ref{Ome}) makes sense also for general non-polynomial Hamiltonians. Consequently, according to \cite{BIT10}, we give the following definition. 

\begin{defi} The Runge-Kutta method  (\ref{RK})--(\ref{Ome})  is called HBVM$(k,s)$.\end{defi} 
The following properties \cite{BIT10,BIT12_1} elucidate the role of two indices $k$ (number of ascissae) and $s$ (degree of the underlying polynomial $\sigma$) in the previous definition.

\begin{theo}\label{propHBVM}
For all $k\ge s$, a HBVM$(k,s)$ method:
\begin{itemize}
\item has order $2s$, that is:
$$y_1-y(h) = O(h^{2s+1});$$
\item is energy conserving for all polynomial Hamiltonians of degree $\nu$ satisfying (\ref{nu});
\item for general  non-polynomial (but suitably regular) Hamiltonians, one has:
\begin{equation}\label{order2k}
H(y_1) - H(y_0) = O(h^{2k+1}).
\end{equation}
\end{itemize}
\end{theo}

\begin{rem} From (\ref{order2k}) one deduces that a HBVM$(k,s)$ method is {\em practically} energy-conserving also in the case of non-polynomial Hamiltonians, provided that $k$ is large enough. Indeed, on a computer, it is enough to approximate the involved integrals to within round-off errors.
\end{rem}
\begin{rem} Though the method (\ref{RK})--(\ref{Ome}) has been derived in the context of Hamiltonian systems, we stress that it makes sense also when replacing problem (\ref{ivp}) by a generic (i.e., non Hamiltonian) initial value problem in the form $y' = f(t,y)$ \cite{BIT12_1}. 
\end{rem}

For sake of completeness, and for later reference, we also  report the following result, which actually shows that HBVM$(k,s)$ methods, with the choice (\ref{ci}) of the abscissae, can be regarded as a generalization of the $s$-stage Gauss-Legendre collocation formulae \cite{BIT10}.

\begin{theo}\label{Gauss_s} HBVM$(s,s)$ coincides with the $s$-stage Gauss-Legendre collocation method.\end{theo}

If we set $\bfy$ the (block) vector with the internal stages (\ref{Yi}) and $\bfe=(1,\dots,1)^T\in\RR^k$, the discrete problem generated by a HBVM$(k,s)$ method is given by
\begin{equation}\label{RKprob}
\bfy = \bfe\otimes y_0 + hA\otimes J\, \nabla H(\bfy),
\end{equation}
which is a nonlinear system of (block) dimension $k$.\footnote{Here $\nabla H(\bfy)$ is the block vector whose entries are given by the gradient of $H$ evaluated at the $k$ stages.}  However, in view of (\ref{order2k}), $k$ is likely to be much larger than $s$ and, consequently, such a formulation is in general not recommendable. 

To derive a more efficient formulation, let us set $\bfgamma$ the (block) vector containing the coefficients defining the polynomial $\sigma$ in (\ref{sigma}), thus obtaining:
$$
\bfgamma = \PP_s^T\Omega\otimes J\, \nabla H(\bfy),\qquad
\bfy           = \bfe\otimes y_0 + h\PP_{s+1}\XX_s\otimes I\,\bfgamma.
$$
Combined together, such equations provide us with the following discrete problem,
\begin{equation}\label{sizes}
F(\bfgamma) \,\equiv\,\bfgamma - \PP_s^T\Omega\otimes J \, \nabla H\left(\bfe\otimes y_0 + h\PP_{s+1}\XX_s\otimes I\,\bfgamma\right) = \bf0,
\end{equation}
whose (block) size is always $s$, independently of $k$.   In general, quite inexpensive  iterations (e.g., the fixed-point iteration) could be used for solving (\ref{sizes}). Nevertheless, in case, e.g., of {\em stiff oscillatory problems}, this could not be practical, since a very small stepsize $h$ would be required: in such a case, a Newton-type iteration is more appropriate  (see the second test problem in Section~\ref{five}). As a popular example,  one easily checks that the simplified Newton iteration, applied for solving (\ref{sizes}), consists in the following iteration \cite{BIT11}:
\begin{eqnarray}\label{simpNewt}
{\rm solve:~ } \left[  I - h X_s\otimes J \nabla^2 H_0\right] \Delta^\ell &=& -F(\bfgamma^\ell)\\
\bfgamma^{\ell+1} &=& \bfgamma^\ell + \Delta^\ell, \qquad \ell=0,1,\dots,\nonumber
\end{eqnarray}
where $\nabla^2 H_0$ is the Hessian of $H(y)$ evaluated at $y_0$. Consequently, the bulk of the computational cost is due to the factorization of the matrix
$$M_0 = I - h X_s\otimes J \nabla^2 H_0,$$ having dimension $2sm\times 2sm$. In the next section, we shall see how to efficiently solve the iteration (\ref{simpNewt}).

\begin{rem}\label{GaussVsHBVM}
As is clear from the previous arguments, HBVM$(k,s)$ methods, with a suitable choice of $k$, is (at least practically) energy-conserving, whereas HBVM$(s,s)$ (i.e., the symplectic $s$-stage Gauss method)  in general is not. Consequently, by taking into account that the (block) dimension of the discrete problem generated by a HBVM$(k,s)$ method is always $s$ independently of $k$, this method is preferable to the $s$-stage Gauss method when an accurate conservation of the energy is required.
\end{rem}

\section{The new splitting procedure}\label{three}
The iteration (\ref{simpNewt}) is similar in structure to  the simplified-Newton iteration applied to the original system (\ref{RKprob}), for which a number of splitting procedures have been devised:  as an example, triangular splittings are defined in \cite{HoSw97,HoSw97a,AmBr97,BIM12}; a diagonal splitting, derived from the so called {\em blended implementation} of the methods, is studied in \cite{BrMa02,BrMa04};  additional approaches are described, e.g., in \cite{Bi77,CoVi93,GPGCMo94,HM98,HS98,SWG06}; moreover, we mention that a comprehensive linear analysis of convergence for such iterations (generalizing that at first proposed in \cite{HoSw97}) is reported in \cite{BrMa09}. However, the triangular splitting iteration defined in \cite{HoSw97,HoSw97a}, along with the modified  triangular splitting iteration defined in \cite{AmBr97}, turn out to be not effective for (\ref{simpNewt}), due to the particular structure of the matrix $X_s$ (see (\ref{Xs})). Conversely, the {\em blended iteration} defined in \cite{BrMa02,BrMa04} (see also \cite{BrMa09}), turns out more appropriate, as is shown in \cite{BIT11}. We here shall devise a different iterative procedure, which appears to be even more favourable. This is the subject of the remaining part of this section. The main idea is similar to that explained in \cite{BIM12} for Radau IIA collocation methods, even though the framework and the overall details (and results) are definitely different:  i.e., to replace the set of $s$ (block) unknowns, given by entries of the (block)  vector $\bfgamma$ defined in (\ref{sizes}), with a more convenient one.
To begin with, let us consider the polynomial (\ref{sigma1})
and introduce the new set of (block) unknowns, 
\begin{equation}
\label{hgi}
\hg_i \equiv \sum_{j=0}^{s-1} P_j(\hc_i) \gamma_j(\sigma), \qquad i=1,\dots,s,
\end{equation}
defined as the evaluation of (\ref{sigma1}) 
 at the set of  distinct {\em auxiliary abscissae}
\begin{equation}\label{hci}
\hc_1\,,\,\dots\,,\,\hc_s.
\end{equation}
Introducing the (block) vector
\begin{equation}\label{hbg}
\hbg = \pmatrix{c} \hg_1\\ \vdots\\ \hg_s\endpmatrix,
\end{equation} and the matrix 
\begin{equation}\label{hP}
\hPP = \pmatrix{c} P_{j-1}(\hc_i)\endpmatrix\in\RR^{s\times s},
\end{equation} 
we can recast (\ref{hgi}) in vector form as 
\begin{equation}\label{hgg}
\hbg = \hPP\otimes I\, \bfgamma.
\end{equation}
In terms of the new unknown vector $\hbg$, the simplified Newton iteration (\ref{simpNewt}) reads:
\begin{eqnarray}\label{simpNewt1}
{\rm solve:~ } \hM_0 \hDelta^\ell &=& -\hPP\otimes I\, F(\hPP^{-1}\otimes I\,\hbg^\ell) \equiv \bfeta^\ell,\\[2mm]
\hbg^{\ell+1} &=& \hbg^\ell + \hDelta^\ell, \qquad \ell=0,1,\dots,\nonumber
\end{eqnarray}
where
\begin{equation}\label{hA}
\hM_0 = I - h \left(\hPP X_s\hPP^{-1}\right)\otimes J \nabla^2 H_0 \equiv   I - h \hA \otimes J \nabla^2 H_0.
\end{equation}
\begin{rem}
We stress that matrix $\hA=\hPP X_s\hPP^{-1}$ is independent of $k$: it only depends on $s$, whichever is the considered value of $k\ge s$. Consequently, the approach presented below also applies to the case $k=s$, that is, to the $s$-stages Gauss method.
\end{rem}
The key idea is that of choosing the abscissae (\ref{hci}) such  that  $\hA$ can be factored as
\begin{equation}\label{crout}
\hA = \hL\hU,
\end{equation}
with $\hU$ upper triangular with unit diagonal entries, and $\hL$ lower triangular with constant diagonal entries. In such a case, by following the approach of van der Houwen et al. \cite{HoSw97,HoSw97a}, the iteration (\ref{simpNewt1}) is replaced by the {\em inner-outer} iteration
\begin{eqnarray}\nonumber
{\rm solve:~ } \left[I-h\hL\otimes J\nabla^2H_0\right] \hDelta^{\ell,r+1} &=& h\hL(\hU-I)\otimes J\nabla^2H_0\,\hDelta^{\ell,r}+ \bfeta^\ell, \\ &&r=0,1,\dots,\mu-1, \label{simpNewt2}\\[1em]
\hbg^{\ell+1} &=& \hbg^\ell + \hDelta^{\ell,\mu}, \qquad \ell=0,1,\dots.\nonumber
\end{eqnarray}
In particular, since $\hat \Delta^{\ell,0}=0$, the choice $\mu = 1$ corresponds to the approach used by van der Houwen et al.  to devise PTIRK methods \cite{HoSw97}, whereas, if $\mu$ is large enough to have full convergence of the inner-iteration (the one on $r$), then the outer iteration is equivalent to (\ref{simpNewt1}). Clearly, all the intermediate possibilities can be suitably considered.  

After the convergence of (\ref{simpNewt2}), the new approximation is computed (see (\ref{sigma})) as $$y_1=y_0+h\gamma_0,$$ where $\gamma_0$  (i.e., the first block entry of the vector $\bfgamma$), is retrieved from (\ref{hgg}).  We observe that the diagonal entries of the factor $\hL$ are all equal to a given value, say $d_{s}$, has the obvious advantage that one only needs to factor the matrix
\begin{equation}\label{Omega0}
I-h d_s J\nabla^2 H_0\in\RR^{2m\times 2m}.
\end{equation}

\begin{rem}
 We observe that, in an actual computational code, such a matrix can be kept constant over a number of steps, being factored only when the Hessian needs to be revaluated and/or the stepsize is modified. In this paper, we deliberately ignore this issue, which requires a further analysis (see, e.g., \cite{BrMa05} for the code described in \cite{BrMa04}). Consequenlty, in the numerical tests we shall use a constant stepsize and compute the Hessian at each step.\end{rem}

Concerning $d_s$, the following result holds true.

\begin{theo} Assume that the factorization (\ref{crout}) is defined and that the factor $\hL$ has all its diagonal entries equal to $d_s$. Then, with reference to (\ref{xi}), one has:
\begin{equation}\label{ds}
d_s = \left\{\begin{array}{ccc} \,^s\sqrt{ \prod_{i=1}^{\lfloor \frac{s}2\rfloor} \xi_{2i-1}^2}\,,&~& \mbox{if $s$ is even,}
\\[2mm]
 \,^s\sqrt{\frac{1}2 \prod_{i=1}^{\lfloor \frac{s}2\rfloor} \xi_{2i}^2}\,,&~& \mbox{if $s$ is odd.}\end{array}\right.
 \end{equation}
 \end{theo}
 
 \proof
 Assume that (\ref{hA})--(\ref{crout}) hold true. Then
$$
 \det(X_s) = \det( \hPP X_s\hPP^{-1} ) = \det(\hA) = \det(\hL\hU) = \det(\hL) = d_s^s,
$$
 since $\hU$ has unit diagonal and all the entries of $\hL$ are equal to $d_s$. Consequently,
 $$d_s = \,^s\sqrt{\det(X_s)}.$$ The thesis then follows by considering that, from (\ref{Xs}), $$\det(X_1) = \frac{1}2, \qquad \det(X_2) = \xi_1^2,$$ and, by applying the Laplace expansion, one obtains:
\begin{equation}\label{detXs}
\det(X_s) = \left\{\begin{array}{ccc}  \prod_{i=1}^{\lfloor \frac{s}2\rfloor} \xi_{2i-1}^2\,,&\quad& \mbox{if $s$ is even,}
\\[2mm]
 \frac{1}2 \prod_{i=1}^{\lfloor \frac{s}2\rfloor} \xi_{2i}^2\,,&\quad& \mbox{if $s$ is odd.}\end{array}\right.
\end{equation}
\QED

\medskip
By virtue of the previous result, in order to compute the auxiliary abscissae (\ref{hci}), we have symbolically solved the following set of equations, which is obviously equivalent to requiring that the factor $\hL$ has the diagonal entries equal to each other:
\begin{equation}\label{conds}
\det(\hA_{\ell+1}) = d_s\det(\hA_{\ell}), \qquad \ell=1,\dots,s-1,
\end{equation}
 where $\hA_\ell$ denotes the principal leading submatrix of order $\ell$ of $\hA$, and $d_s$ is given by (\ref{ds}).
\begin{rem}\label{hcs}
We observe that the auxiliary abscissae (\ref{hci}) are $s$, whereas the  algebraic conditions (\ref{conds}) are $s-1$.  This means that a further condition can be imposed on the abscissae: we have chosen it in order to  improve the convergence properties of the iteration (\ref{simpNewt2}), according to the linear analysis of convergence reported in Section~\ref{four}; in particular, we shall
(approximately) minimize the {\em maximum amplification factor} of the iteration,  as it will be later explained.
\end{rem}
The obtained results  are listed in Table~\ref{cds}, for $s=2,\dots,6$, from which one sees that in all cases the abscissae are distinct and inside the interval $[0,1]$.

We emphasize that, for any given $s$, the distribution of the auxiliary abscissae (\ref{hci}) 
is independent of $k$ and so is the factorization (\ref{crout}) of the matrix $\hA$ whose computation is responsible of the bulk of the computational effort during the integration process.  This property has a relevant consequence during the implementation phase of this class of methods. In fact, one can  conjecture a procedure to advance the time that dynamically selects the most appropriate value of $k$. Depending on the specific problem at hand and the configuration of the system at the given time, one can easily switch from a symplectic to an energy preserving method by choosing $k=s$ (Gauss method) or $k>s$, respectively.          
\begin{table}[p]
\caption{Auxiliary abscissae (\ref{hci}) for the HBVM$(k,s)$ and $s$-stage Gauss method, $s=2,\dots,6$,
and the diagonal entry $d_{s}$ (see (\ref{ds})) of the corresponding factor $\hL$.}
\label{cds}
\vspace{2mm}
\centerline{\begin{tabular}{|c|c|}
\hline
\multicolumn{2}{|c|}{$s=2$}\\
\hline
$\hc_1$ & 0.26036297108184508789101036587842555\\
$\hc_2$ & 1\\
\hline
$d_2$    & 0.28867513459481288225457439025097873\\
\hline
\hline
\hline
\multicolumn{2}{|c|}{$s=3$}\\
\hline
$\hc_1$ & 0.15636399930006671060146617869938122 \\
$\hc_2$ & 0.45431868644630821020177903150137523 \\
$\hc_3$ & 0.948\\
\hline
$d_3$    & 0.20274006651911333949661483325792675 \\
\hline
\hline
\multicolumn{2}{|c|}{$s=4$}\\
\hline
$\hc_1$ & 0.11004843257056123468614502691988075 \\
$\hc_2$ & 0.31588689139705398683980065724981436 \\
$\hc_3$ & 0.53114668286639796587351917750274705 \\
$\hc_4$ & 0.884 \\
\hline
$d_4$   & 0.15619699684601279005430416526875577 \\
\hline
\hline
\multicolumn{2}{|c|}{$s=5$}\\
\hline
$\hc_1$ & 0.084221784434612320884185541600934218\\ 
$\hc_2$ & 0.248618520588562018051811779022293944\\
$\hc_3$ & 0.413725268815220956415498643302145284\\
$\hc_4$ & 0.587098748971877116030882436751962384\\
$\hc_5$ & 0.9338 \\
\hline
$d_5$    & 0.12702337351164258963093490787943281 \\
\hline
\hline
\multicolumn{2}{|c|}{$s=6$}\\
\hline
$\hc_1$ & 0.20985774196263657630356114041757724\\ 
$\hc_2$ & 0.36816786358152563671526302698797908\\
$\hc_3$ & 0.39607328223635472401921951140390213\\
$\hc_4$ & 0.62783521091780460858476326939502046\\
$\hc_5$ & 0.04580307227138364391540767310611717\\
$\hc_6$ & 0.94225 \\
\hline
$d_6$    & 0.10702845478806509529222890981996019 \\
\hline
\end{tabular}}
\end{table}

\section{Convergence analysis and comparisons}\label{four}

In this section we briefly analyze the splitting procedure (\ref{simpNewt2}).  In general, its convergence properties could be discussed in the framework of quasi-Newton methods, leading to the (quite) obvious result that linear convergence is obtained for sufficiently small $h$, provided that $H$ is suitably regular. Nevertheless, a more suited approach, which has proved to be very effective in the actual design of efficient variable-order/variable-stepsize codes for ODE-IVPs (see, e.g., \cite{BrMa04}), is based on the linear analysis of convergence in \cite{HoSw97} (further developed in \cite{BrMa09}). Such an analysis is well motivated from the fact that the inner iteration in (\ref{simpNewt2}) amounts to solving a linear system. This latter system can be thought of being obtained by applying the original numerical method to the local (frozen) linearized problem. As a consequence, one can decompose it in the subspaces spanned by  the eigenvalues of the Jacobian. Equivalently, one can directly consider the scalar problem defined by each eigenvalue. Consequently,  one is led to study the behavior of the method when applied to the celebrated test equation:
\begin{equation}\label{test}
y' = \lambda y, \qquad y(t_0)=y_0.
\end{equation}
Clearly, one directly arrives to the same conclusion in case problem (\ref{ivp}) is separable, with a quadratic Hamiltonian, with the eigenvalues lying on the imaginary axis. Since problem (\ref{test}) is linear, the iteration (\ref{simpNewt2}) consists in solving the inner iteration alone, so that we can skip the index $\ell$ of the outer iteration. By setting, as  is usual, $q=h\lambda$, one then obtains that the error equation associated with the iteration (\ref{simpNewt2}) is given by
\begin{equation}\label{erreq}
e_{r+1} = Z(q)e_r, \qquad Z(q):= q(I-q\hat{L})^{-1}\hat{L}(\hat{U}-I), \qquad r=0,1,\dots,
\end{equation}
where  $e_r$ is the error vector at step $r$ and $Z(q)$ is the iteration matrix induced by the splitting procedure. This latter will converge if and only if its spectral radius,
$$\rho(q) := \rho(Z(q)),$$
is less than $1$. The {\em region of convergence} of the iteration is then defined as
$$\DD = \left\{ q\in\CC\,:\, \rho(q)<1\right\}.$$
The iteration is said to be {\em $A$-convergent} if $\CC^-\subseteq\DD$. If, in addition, the {\em stiff amplification factor},
$$\rho^\infty := \lim_{q\rightarrow\infty} \rho(q),$$
is null, then the iteration is said to be {\em $L$-convergent}.\,\footnote{In general, $A$-convergent iterations are appropriate for $A$-stable methods, and $L$-convergent iterations are appropriate for $L$-stable methods.}
In our case, since
\begin{equation}\label{Minfty}
Z(q) \rightarrow (\hat{U}-I), \qquad q\rightarrow\infty,
\end{equation}
which is a nilpotent matrix of index $s$, the iteration is $L$-convergent if and only if it is $A$-convergent.
Since the iteration is well defined for all $q\in\CC^-$ (due to the fact that the diagonal entry of $\hat{L}$, $d_s$, is positive, as was shown in (\ref{ds}))
 and $\rho(0)=0,$ from the maximum-modulus theorem it follows immediately that $A$-convergence is, in turn,  equivalent to require that the {\em maximum amplification factor},
$$\rho^* := \max_{x\in\RR} \rho(ix),$$
is not larger than $1$. Another useful parameter is the {\em nonstiff amplification factor},
\begin{equation}\label{M0}
\tilde\rho := \rho(  \hat{L}(\hat{U}-I) ),\end{equation}
that governs the convergence of the iteration for small values of $q$, since
$$\rho(q) \approx \tilde\rho |q|, \qquad \mbox{for} \quad q\approx 0.$$
Clearly, the smaller $\rho^*$ and $\tilde\rho$, the better the convergence properties of the iteration.

With these premises, we can now better specify what anticipated in Remark~\ref{hcs}, concerning the additional condition imposed to derive the auxiliary abscissae (\ref{hci}), while fulfilling the conditions (\ref{conds}). In more details, the abscissae listed in Table~\ref{cds} have been computed by (approximately) solving the following constrained minimization problem:
\begin{eqnarray*}
\min_{\hc_1,\dots,\hc_s} && \rho^*\\
s.t. && \det(\hA_{\ell+1}) = d_s\det(\hA_\ell), \qquad \ell=1,\dots,s-1.
\end{eqnarray*}
Clearly, this has been made possible by the introduction of the transformation (\ref{hgg}).

In Table~\ref{parametri} we list the maximum amplification factors and the nonstiff amplification factors
 for the following $L$-convergent iterations applied to the $s$-stage Gauss-Legendre methods:
\begin{itemize}
\item[(i)] the iteration obtained by the original triangular splitting in \cite{HoSw97};
\item[(ii)] the iteration obtained by the modified triangular splitting in \cite{AmBr97};
\item[(iii)] the {\em blended} iteration obtained by the {\em blended implementation} of the methods, as defined in \cite{BrMa02};
\item[(iv)] the iteration defined by (\ref{simpNewt2}).

\end{itemize}
We recall that the scheme (i) (first column) requires $s$ real factorizations per iteration, whereas (ii)--(iv) only need one factorization per iteration.  From the parameters listed in the table, one concludes that the proposed splitting procedure is the most effective among all the considered ones.

\begin{rem}
For sake of accuracy, we stress that, when dealing with the actual implementation of HBVM$(k,s)$ methods, only the blended iteration \cite{BIT11}  and the newly proposed one (\ref{simpNewt2}) can be considered,  whereas the triangular splitting defined in \cite{HoSw97} and its modified version \cite{AmBr97} turn out to be not effective, as was pointed out at the beginning of Section~\ref{three}. Consequently, in such a case, one has to consider only the last two group of columns in Table~\ref{parametri}.
\end{rem}

\begin{table}[t]
\begin{center}
\caption{Amplification factors for the triangular splitting in \cite{HoSw97},
the modified triangular splitting in \cite{AmBr97}, the {\em blended} iteration in \cite{BIT11},
 and the splitting (\ref{simpNewt2}), for the $s$-stage Gauss-Legendre formulae. The last two cases coincide with those for the HBVM$(k,s)$ methods, $k\ge s$.}
\label{parametri}

\vspace{2mm}
\begin{tabular}{|c|cc|cc|cc|cc|}
\hline
       &\multicolumn{2}{c|}{(i): triangular}
       &\multicolumn{2}{c|}{(ii): triangular}
       &\multicolumn{2}{c|}{(iii): {\em blended}}
       &\multicolumn{2}{c|}{(iv): triangular}\\
       &\multicolumn{2}{c|}{splitting in \cite{HoSw97}}
       &\multicolumn{2}{c|}{splitting in \cite{AmBr97}}
       &\multicolumn{2}{c|}{iteration in \cite{BIT11}}
       &\multicolumn{2}{c|}{splitting  (\ref{simpNewt2})}\\
\hline 
$s$ &  $\rho^*$ & $\tilde{\rho}$ & $\rho^*$  & $\tilde{\rho}$ & $\rho^*$ & $\tilde{\rho}$ & $\rho^*$ & $\tilde{\rho}$ \\
\hline
2   & 0.1429 & 0.0833 & 0.1340 & 0.0774 & 0.1340 & 0.0774 & 0.1340 &   0.0774 \\
3  &  0.3032 & 0.1098 & 0.2537 & 0.0856 & 0.2765 & 0.1088 & 0.2536 &   0.0870 \\
4  &  0.4351 & 0.1126 & 0.3492 & 0.0803 & 0.3793 & 0.1119  & 0.3291 &   0.0859 \\
5  &  0.5457 & 0.1058 & 0.4223 & 0.0730 & 0.4544 & 0.1066 &  0.3709 &   0.0654 \\
6  &  0.6432 & 0.0973 & 0.4861 & 0.0702 & 0.5114 & 0.0993 &  0.4353 &    0.0650 \\
\hline
\end{tabular}\end{center}
\end{table}

\subsection{Averaged amplification factors}
The previous amplification factors measure the asymptotic speed of convergence when an infinite number of iterations are performed, which is not the case, in the actual implementation of the methods. For this purpose (see, e.g., \cite{BrMa09}) it is also customary to define corresponding {\em averaged} amplification factors, which measure the ``average'' convergence when a prescribed number of iterations is performed. In particular, by considering a suitable matrix norm $\|\cdot\|$, and with reference to what previously has been set out, we define the following {\em averaged amplification factors} when $\mu$ iterations of (\ref{erreq}) are carried out:
\begin{equation}\label{averfat}
\rho_\mu^*:=\sup_{x\in\RR}\,^\mu\sqrt{\|Z(ix)^\mu\|},\qquad \tilde{\rho}_\mu:=\,^\mu\sqrt{\left\|\left[\hat{L}(\hat{U}-I)\right]^\mu\right\|},\qquad \rho_\mu^\infty:=\,^\mu\sqrt{\|(\hat{U}-I)^\mu\|}.
\end{equation}
Clearly, 
$$\lim_{\mu\rightarrow\infty} \rho_\mu^* = \rho^*, \qquad \lim_{\mu\rightarrow\infty} \tilde\rho_\mu = \tilde\rho,$$
and $$\rho_\mu^\infty = 0, \qquad \forall \mu\ge s.$$
In Table~\ref{averpar} we list the obtained averaged amplification factors (\ref{averfat}) when performing $\mu=1,2,3$ iterations, by considering the infinity norm. As one may see, the resulting iteration turns out to be $A$-convergent also when using just one inner iteration, unless the case $s=6$, which requires at least 3 inner iterations.

\begin{table}[t]
\begin{center}
\caption{Averaged amplification factors (\ref{averfat}) for the splitting (\ref{simpNewt2}), used for the HBVM$(k,s)$ methods, $k\ge s$, when performing $\mu=1,2,3$ iterations.}
\label{averpar}

\vspace{2mm}
\begin{tabular}{|c|ccc|ccc|ccc|}
\hline 
$s$ &  $\rho_1^*$ & $\tilde\rho_1$ & $\rho_1^\infty$  & $\rho_2^*$ & $\tilde\rho_2$ & $\rho_2^\infty$ & $\rho_3^*$ & $\tilde\rho_3$ & $\rho_3^\infty$ \\
\hline
2   &  0.1340  &  0.0774 &  0.0981 &  0.1340  & 0.0774  &  0 &   0.1340 &  0.0774 &  0  \\
3  &  0.4492  &  0.0874  &  0.2606 &   0.3423 &  0.0873 &  0.1091 &  0.3087 &  0.0872 &     0\\
4  &   0.4751  &  0.1459 &   0.4751 &   0.4098 &   0.1200 &   0.1757 &   0.3848 &   0.1091 &   0.1294\\
5  &  0.8625   & 0.2045  &  0.7471  &  0.6775  &  0.1385  &  0.2872  &  0.5874  &  0.1154  &  0.1747\\
6  &  3.0797   & 0.2747  &  1.4988  &  1.2780  &  0.1356  &  0.4929  &  0.9451  &  0.1121  &  0.2697\\
\hline
\end{tabular}\end{center}
\end{table}

\begin{rem} When performing only $\mu$ inner-iterations for solving the discrete problem generated by (\ref{test}), we have to consider also the {\em outer} iteration, even though the problem is linear. In such a case, by setting $E_\ell$ the error at the $\ell$-th outer iteration, it is quite straightforward to see that the error equation for the outer iteration is given by (compare with (\ref{simpNewt2})):
$$E_{\ell+1} = Z(q)^\mu E_\ell, \qquad \ell=0,1,\dots.$$
Consequently, the previous convergence analysis also applies to the present case.
\end{rem}

\section{Numerical Tests}\label{five}
In this section, we report a couple of numerical examples, aimed to put into evidence the features of the methods, and/or the effectiveness of the proposed iterative procedure.  For both problems, we list the computational cost for HBVM$(k,s)$ methods, in terms of required iterations for solving the generated discrete problems with a constant stepsize, when using:
\begin{itemize}
\item the fixed-point iteration;
\item the blended iteration in \cite{BIT11};
\item the splitting iteration (\ref{simpNewt2}) with 2 inner iterations.
\end{itemize}
The choice of $2$ inner iterations in  (\ref{simpNewt2}) makes the cost of one outer iteration comparable to that of one blended iteration, provided that (\ref{simpNewt2}) is implemented as suggested in \cite{BIM12}. The total number of functional evaluations equals the number of iterations times $k$. Moreover, for the latter two iterations, at each step one also needs to evaluate the Hessian $\nabla^2H$, as well as to factor a matrix having the same size as that of the continuous problem (i.e., 
(\ref{Omega0}), in the case of the iteration (\ref{simpNewt2})).

The first problem is a nonlinear Hamiltonian problem describing the motion of a charged particle, with charge $e$ and mass $m$, in a magnetic field with Biot-Savart potential. 
It is defined by the Hamiltonian:
\begin{equation}\label{biot}
H(x,y,z,x',y',z')= \frac{1}{2m}\left[ \left( x'-\aa\frac{x}{\rho^2}\right)^2 +
\left( y'-\aa\frac{y}{\rho^2}\right)^2 +\left(z'+\aa\log\rho\right)^2
\right],
\end{equation}
with $\rho=\sqrt{x^2+y^2}$ and $\aa=e B_0$, $B_0$ being the intensity of the magnetic field. We have used the values
$$m=1, \qquad e=-1, \qquad B_0=1,$$ and the initial values
\begin{equation}\label{biot0}
x=0.5, \qquad y=10, \qquad x'=-0.1, \qquad y'=-0.3, \qquad z=z'=0.
\end{equation}
In Table~\ref{biot_tab} we list the results obtained by applying the HBVM$(k,2)$ methods, $k=2,4,6,8,10$, for solving this problem over the interval $[0,10^3]$ with stepsize $h=0.1$. From the results in the table, one infers that:
\begin{itemize}
\item the Hamiltonian error monotonically decreases as $k$ is increased and, for $k=10$, one obtains a practical conservation, for the given stepsize (consequently, larger values of $k$ would be useless);

\item the solution error when using the symplectic 2-stages Gauss method (i.e., HBVM(2,2)) is larger than that obtained when the energy error decreases;

\item the proposed iterative procedure (\ref{simpNewt2}) is more effective than the blended iteration proposed in \cite{BIT11}. In such a case, however, both iterations turn out to be not very competitive, with respect to the use of a fixed-point iteration, since this problem is not {\em stiff};

\item all iterations provide a total cost which is independent of $k$.
\end{itemize}
\begin{table}[t]
\caption{Results when solving Problem (\ref{biot})-(\ref{biot0}) by using the HBVM$(k,2)$ method with stepsize $h=0.1$ over the interval $[0,10^3]$.}
\label{biot_tab}
\centerline{\begin{tabular}{|r|l|l|r|r|r|}
\hline
       & Hamiltonian & solution & fixed-point & blended   & splitting \\
$k$ & error            & error      & iterations   & iterations & iterations\\
\hline
2 &  $1.6\cdot10^{-3}$ & $9.97\cdot10^{-2}$ &  79511   & 66854 &  48030\\
4 &  $8.3\cdot10^{-6}$ & $1.82\cdot10^{-2}$ &  79846   & 66884 & 48252\\
6 &  $5.9\cdot10^{-9}$ & $1.81\cdot10^{-2}$  & 79911   & 66941 & 48349\\
8 &  $1.7\cdot10^{-12}$ & $1.81\cdot10^{-2}$ & 79939  & 66963 & 48377\\
10 & $4.4\cdot10^{-16}$ & $1.81\cdot10^{-2}$ & 79962 & 66976 & 48402\\
\hline
\end{tabular}}
\end{table}

The second test problem that we consider is, on the contrary, a {\em stiff oscillatory} problem. It is defined as a slight modification of the Fermi-Pasta-Ulam problem described in \cite{HLW06}.\footnote{The original problem reported in \cite{HLW06} is obtained by  setting $m=3$ and $\omega_i=50$, $i=1,\dots,m$, in (\ref{fpumod}).} The Hamiltonian is now given by:
\begin{equation}\label{fpumod}
H(p,q) = \frac{1}2\sum_{i=1}^m\left(p_{2i-1}^2+p_{2i}^2\right)
+\frac{1}4\sum_{i=1}^m\omega_i^2\left(q_{2i}-q_{2i-1}\right)^2
+\sum_{i=0}^m\left(q_{2i+1}-q_{2i}\right)^4,
\end{equation}
\no with $q,p\in\RR^{2m}$  ~and~ $q_0=q_{2m+1}=0$.  We choose ~$m=7$,~ so that the problem has dimension 28, and
\begin{equation}\label{wi}
\omega_i=\omega_{m-i+1}=10, \qquad i=1,2,3,  \qquad\mbox{and}\qquad \omega_4 = 10^4.
\end{equation}
The starting vector is
\begin{equation}\label{fpu0}
p_i=0, \qquad q_i = \frac{i-1}{2m-1}, \qquad i=1,\dots,2m.
\end{equation}
In such a case, the Hamiltonian function is a polynomial of degree 4, so
that the  HBVM$(2s,s)$ method (having order $2s$), is able to exactly preserve the Hamiltonian.
As an example, fix $s=3$ and integrate the problem on the interval $[0,10]$. In this case, the fixed-point iteration cannot be expected to work, when using stepsizes much larger than $\|\omega\|_\infty^{-1}=10^{-4}$, as is confirmed by the results listed in Table~\ref{fpu_mod_fpi}. 
Similarly, explicit methods, which exist in this specific case since the problem is separable (see \cite[Chapter~8]{SSC94}), suffer from similar restrictions on the stepsize because of stability reasons. In particular, we consider a composition method, having order 6, based on the second order St\"ormer-Verlet method (see \cite[Chapter~II.4]{HLW06} for details), requiring 18 function evaluations per step\,\footnote{Consequently, each step of this composition method has a cost which is comparable to 3 fixed-point iterations for HBVM(6,3).}: the results listed in Table~\ref{compver6} clearly confirm this fact.

\begin{figure}[hp]
\centerline{\includegraphics[width=0.75\textwidth]{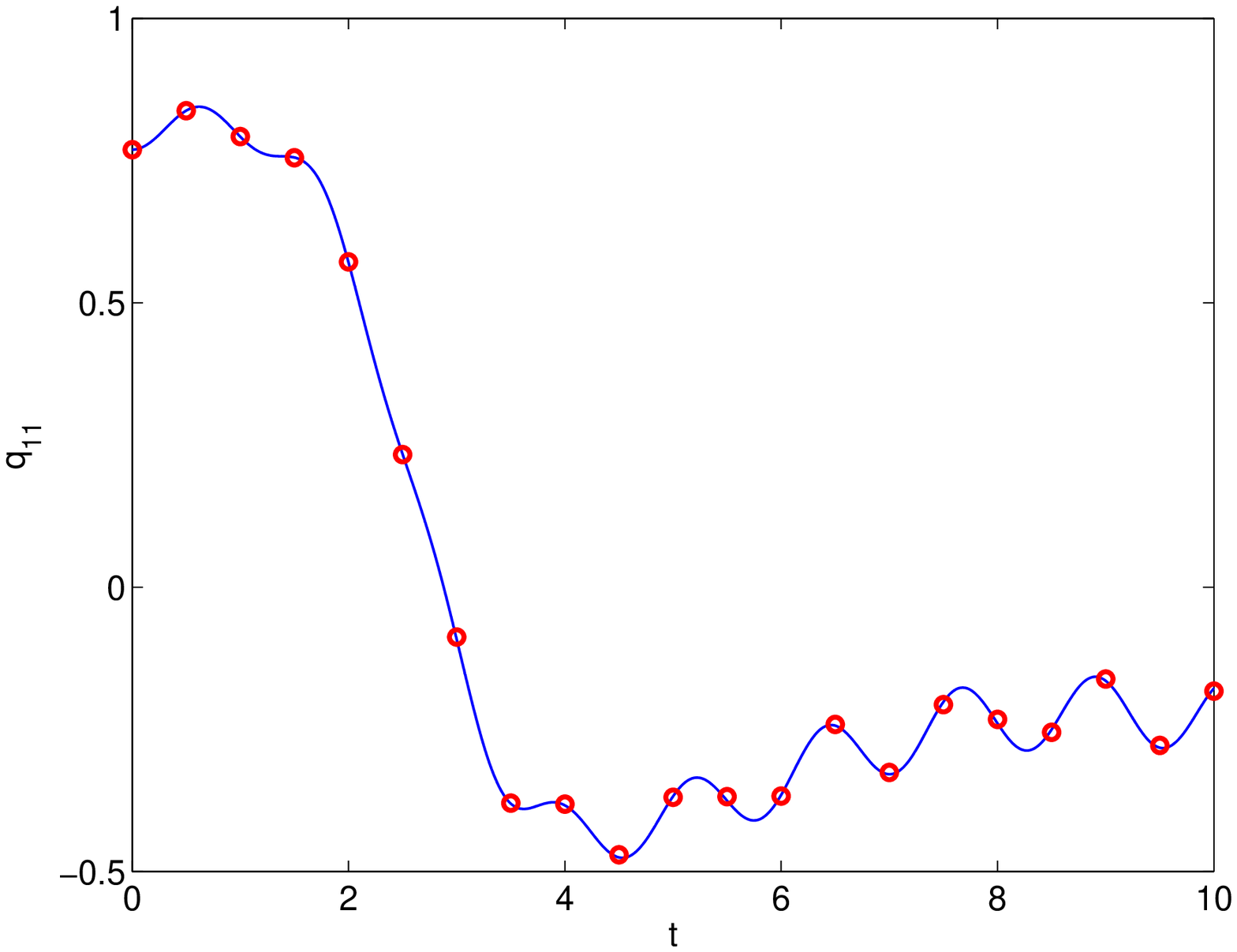}}
\caption{Numerical approximation obtained by using HBVM(6,3) with stepsizes $h=10^{-4}$ (continuous line) and $h=0.5$ (circles) for solving problem (\ref{fpumod})--(\ref{fpu0}).}
\label{q11}
\bigskip

\centerline{\includegraphics[width=0.75\textwidth]{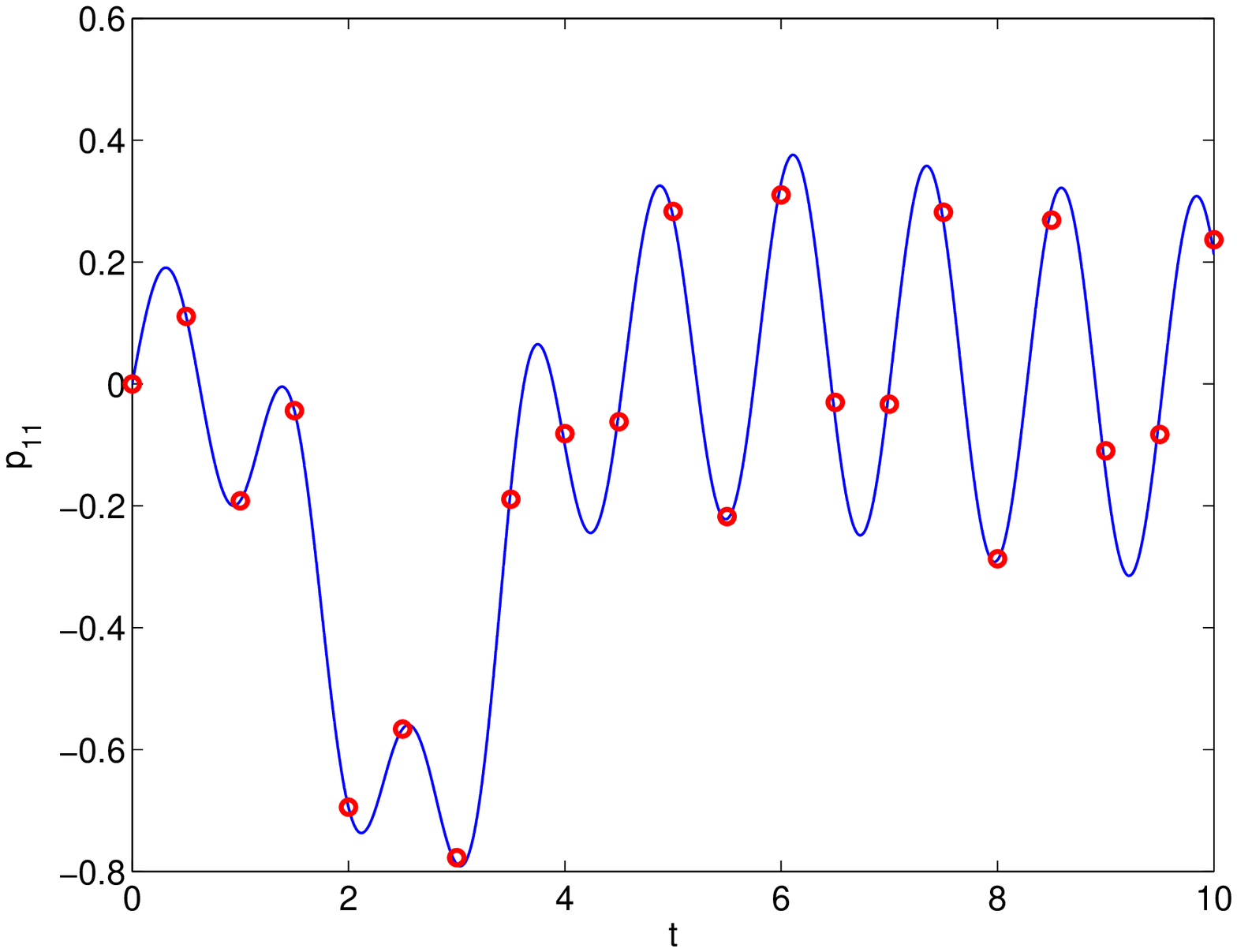}}
\caption{Numerical approximation obtained by using HBVM(6,3) with stepsizes $h=10^{-4}$ (continuous line) and $h=0.5$ (circles)  for solving problem (\ref{fpumod})--(\ref{fpu0}).}
\label{p11}
\end{figure}

Conversely, the use of Newton-type iterations for solving the discrete problems generated by the HBVM(6,3) method, permits to use much larger stepsizes, thus allowing to approximate the low frequencies without being hindered by the high ones. By using the blended iteration defined in \cite{BIT11} and the iteration (\ref{simpNewt2}) previously defined, one obtains the results listed in Table~\ref{fpu_mod_newt}. Even when using very coarse stepsizes, the approximation of the slowly-oscillating components of the solution (24 out of 28) is satisfactory: as an example, in Figures~\ref{q11} and \ref{p11} there is the plot of the slowly-oscillating components $q_{11}$ and $p_{11}$, respectively, by using a finer step, $h=10^{-4}$, and a much coarser one, $h=0.5$.\footnote{By the way, we mention that also the {\em amplitude} of the remaining 4 highly-oscillating components turns out to be well approximated, when using a stepsize $h=0.1$.} Last but not least, from the figures in Table~\ref{fpu_mod_newt}, one sees that the new iterative procedure (\ref{simpNewt2}) is the most effective one, though using only 2 inner iterations.

\begin{table}
\caption{Fixed-point iterations for solving problem (\ref{fpumod})--(\ref{fpu0}), on the interval [0,10], by using HBVM(6,3) with stepsize $h$ (*** means that the iteration doesn't converge).}
\label{fpu_mod_fpi}

\smallskip
\centerline{\begin{tabular}{|r|r|}
\hline
       &fixed-point \\
$h$ &iterations\\
\hline
$10^{-4}$            &   2278912   \\
$2\cdot 10^{-4}$ &   1904534   \\
$4\cdot 10^{-4}$ &  4540389    \\
$5\cdot 10^{-4}$ & *** \\
\hline
\end{tabular}}

\bigskip
\caption{Hamiltonian error, obtained by using a sixth-order explicit composition method based on the St\"ormer-Verlet method, for solving problem (\ref{fpumod})--(\ref{fpu0}) on the interval [0,10] by using stepsize $h$ (*** means that the numerical solution diverges).}
\label{compver6}

\smallskip
\centerline{\begin{tabular}{|r|l|}
\hline
       &Hamiltonian \\
$h$ &error\\
\hline
$10^{-5}$            & $9.2\cdot10^{-8}$     \\
$5\cdot10^{-5}$  & $1.5\cdot10^{-3}$     \\
$10^{-4}$            & $8.5\cdot10^{-2}$     \\
$2\cdot 10^{-4}$ & *** \\
$4\cdot 10^{-4}$ & *** \\
$5\cdot 10^{-4}$ & *** \\
\hline
\end{tabular}}

\bigskip

\caption{Newton-type iterations for solving problem (\ref{fpumod})--(\ref{fpu0}), on the interval [0,10], by using HBVM(6,3) with stepsize $h$.}
\label{fpu_mod_newt}
\centerline{\begin{tabular}{|r|r|r|}
\hline
        & blended   & splitting \\
$h$  & iterations & iterations\\
\hline
$10^{-4}$            & 1634792 &  856691  \\
$5\cdot 10^{-4}$ &   599927 &  299586  \\
$10^{-3}$            &  241468 &  141506  \\
$5\cdot 10^{-3}$ &    29051 &    19148  \\
$10^{-2}$            &    12721 &      8955 \\
$5\cdot 10^{-2}$ &      2369 &      1556 \\
$10^{-1}$            &      1400 &        864 \\
$5\cdot 10^{-1}$ &        440 &        258 \\
\hline
\end{tabular}}
\end{table}

\section{Conclusions}\label{six}

In this paper we have defined an efficient iterative procedure for solving the discrete problems generated by the application of HBVM$(k,s)$ methods, a class of energy-conserving methods for polynomial Hamiltonian dynamical systems. The proposed implementation turns out to improve over that proposed in \cite{BIT11}. Moreover, it also applies to  Gauss-Legendre formulae and the resulting linear convergence analysis shows that the proposed iterative procedure is the most effective, among those based on suitable splittings of the corresponding Butcher array of the methods, known from the literature. A few numerical tests confirm the effectiveness of the proposed iteration when numerically solving stiff oscillatory problems. 

\subsection*{Acknowledgements} The authors wish to thank the anonymous referees, for the useful comments and remarks.

\end{document}